\documentclass[12pt,twoside]{article}
\usepackage{amsmath, amssymb, amsthm}
\newtheorem{lemma}{Lemma}

\newtheorem{definition}{Definition}

\font\fourteenb=cmb10 at 14pt
\begin{document}
\vspace*{-1.0cm}\noindent \copyright
 Journal of Technical University at Plovdiv\\[-0.0mm]\
\ Fundamental Sciences and Applications, Vol. 6, 1998\\[-0.0mm]
\textit{Series A-Pure and Applied Mathematics}\\[-0.0mm]
\ Bulgaria, ISSN 1310-8271\\[+1.2cm]
\font\fourteenb=cmb10 at 14pt
\begin{center}

   {\bf \LARGE Note on the growth of Area functions
   \\ \ \\ \large Peyo Stoilov, Roumyana Gesheva}
\end{center}

\

\

\begin{abstract}
Let
$$
A(r,f)=\int\limits_{0}^{r}\int\limits_{0}^{2\pi}|f'(\rho
e^{i\theta})|^2 \rho\,d\rho\,d\theta\;,
$$
where $r=|z|<1$. $A(r,f)\,$ is the area of $\{\;z\,:\,|z|\leq
r\,\}\,$ under $f$. In this note we prove that
$$
A(r,f)\leq \frac{1}2 \Vert f\Vert _{H^{\infty}}
\int\limits_{0}^{2\pi}|f'(r e^{i\theta})|\,d\theta
$$
and based on this result, a new proof of the inequalities of D. J.
Hallenbeck for the Area functions of multipliers of fractional
Cauchy transforms is given.
\end{abstract}

\section{Introduction.}
Let $\Bbb{D}\,$ denote the unit disk in the complex plane and $\Bbb{T}\,$-
the unit circle. Let $M\,$ be the Banach space of all complex-valued Borel
measures on $\Bbb{T}\,$ with the usual variation norm. For $\alpha>0$, let
${\cal F}_{\alpha}$ denote the family of functions $g\,$ for which there exists
$\mu\in M\,$ such that
\begin{equation}
g(z)=\int\limits_{T}\frac{1}{(1-\overline{\xi}z)^{\alpha}}\,d\mu(\xi),
\quad z\in\Bbb{D}.
\end{equation}
We note that ${\cal F}_{\alpha}\,$ is a Banach space with the natural norm
\[
\Vert g\Vert_{{\cal F}_{\alpha}}=\inf\left\{\Vert \mu\Vert :\mu\in M
\mbox{ such that (1) holds}\;\right\}.
\]
\begin{definition}
Suppose that $f\,$ is holomorphic in $\Bbb{D}$. Then $f\,$ is called a
{\bf multiplier} of ${\cal F}_{\alpha}\,$ if $g\in{\cal F}_{\alpha}
\Rightarrow fg\in{\cal F}_{\alpha}$.
\end{definition}
Let $m_{\alpha}\,$ denote the set of multipliers of ${\cal F}_{\alpha}\,$
and
\[
\Vert f\Vert_{m_{\alpha}}=\sup\left\{\Vert fg\Vert_{{\cal F}_{\alpha}}
:\Vert g\Vert_{{\cal F}_{\alpha}}\leq 1\right\}.
\]
For $f\,$ analytic in $\Bbb{D}\,$ define
$$
A(r,f)=\int\limits_{0}^{r}\int\limits_{0}^{2\pi}|f'(\rho e^{i\theta})|^2 \rho\,d\rho\,d\theta
$$
for $r=|z|<1$. $A(r,f)\,$ is the area of $\{\;z\,:\,|z|\leq r\,\}\,$ under $f$.

In [1] D. J. Hallenbeck determined the sharp growth of $A(r,f)\,$ for $f\in m_{\alpha}$.
In this note we give new proof of these results.
\section{Main Lemma.}
\begin{lemma}
If $f\in H^{\infty}$, then
$$
A(r,f)\leq \frac{1}2 \Vert f\Vert _{H^{\infty}}
\int\limits_{0}^{2\pi}|f'(r e^{i\theta})|\,d\theta
$$
\end{lemma}
Proof.
\begin{eqnarray*}
A(r,f)& = & \int\limits_{0}^{r}\int\limits_{0}^{2\pi}|f'(\rho e^{i\theta})|^2 \rho\,d\rho\,d\theta= 2\pi \int\limits_{0}^{r}\left( \frac{1}{2\pi}
\int\limits_{0}^{2\pi}|f'(\rho e^{i\theta})|^2\,d\theta\right)\rho\,d\rho\\
      & = & 2\pi \int\limits_{0}^{r}\left( \sum\limits_{k=0}^{\infty}(k+1)^2|\hat{f}(k+1)|^2
\rho^{2k}\right)\rho\,d\rho\\
      & = & 2\pi\sum\limits_{k=0}^{\infty} (k+1)^2|\hat{f}(k+1)|^2\int\limits_{0}^{r}\rho^{2k+1}d\rho\\
      & = & \pi\sum\limits_{k=0}^{\infty} (k+1)|\hat{f}(k+1)|^2 \,r^{2(k+1)}
\end{eqnarray*}
Applying Parseval's relations
\begin{eqnarray*}
\int\limits_{0}^{2\pi}\overline{f(re^{i\theta})}\,f'(re^{i\theta})re^{i\theta}\,d\theta=
\sum\limits_{k=0}^n (k+1)\hat{f}(k+1)r^{k+1}\int\limits_{0}^{2\pi}e^{-i(k+1)\theta}
f(re^{i\theta})\,d\theta\\
=2\pi\sum\limits_{k=0}^{\infty} (k+1)|\hat{f}(k+1)|^2 \,r^{2(k+1)}
\end{eqnarray*}
we obtain
$$
A(r,f)=\frac{1}2 \int\limits_{0}^{2\pi}\overline{f(re^{i\theta})}\,f'(re^{i\theta})re^{i\theta}\,d\theta\,,
$$
which implies
$$
A(r,f)\leq \frac{1}2 \sup\limits_{\theta} |f(re^{i\theta})|\int\limits_{0}^{2\pi}
|f'(re^{i\theta})|\,d\theta\leq \frac{1}2 \Vert f\Vert _{H^{\infty}}
\int\limits_{0}^{2\pi}|f'(r e^{i\theta})|\,d\theta\,.
$$
\section{Applications of the Main Lemma.}
It was proved by Hallenbeck and Samotij that
$$
\int\limits_{0}^{2\pi}|f'(r e^{i\theta})|\,d\theta\leq c\Vert f\Vert _{m_{\alpha}}\frac
{1}{(1-r)^{\alpha}}\quad, 0<\alpha<1\,;
$$
$$
\int\limits_{0}^{2\pi}|f'(r e^{i\theta})|\,d\theta\leq c\Vert f\Vert _{m_1}\frac{1}{1-r}\,
\frac{1}{\log\frac{1}{1-r}}\quad,
$$
where $c\,$ denotes a universal constant which may change from line to line or even within a line.

The inequalities
$$
 A(r,f)\leq c\Vert f\Vert_{H^{\infty}}\Vert f\Vert _{m_{\alpha}}\frac
{1}{(1-r)^{\alpha}}\quad, 0<\alpha<1\,;
$$
$$
A(r,f)\leq c\Vert f\Vert_{H^{\infty}}\Vert f\Vert _{m_1}\frac{1}{1-r}\,
\frac{1}{\log\frac{1}{1-r}}\quad
$$
follow immediately from the Main Lemma. The same inequalities were
obtained by D. J. Hallenbeck [1] with the help of more complicated
method.
 \

\

\begin{center}

\end{center}

\

\

\noindent
{\small Department of Mathematics\\
        Technical University\\
        25, Tsanko Dijstabanov,\\
        Plovdiv, Bulgaria\\
        e-mail: peyyyo@mail.bg}

\end{document}